\documentclass[11pt]{amsart}
\usepackage{amssymb}
\usepackage{amsmath}
\usepackage[active]{srcltx}
\usepackage{t1enc}
\usepackage[latin2]{inputenc}
\usepackage{verbatim}
\usepackage{amsmath,amsfonts,amssymb,amsthm}
\usepackage[mathcal]{eucal}
\usepackage{enumerate}
\usepackage[centertags]{amsmath}
\usepackage{graphics}

\newtheorem{theorem}{Theorem}

\newtheorem{corollary}{Corollary}

\begin{document}
\author{George Tephnadze}
\title[partial sums ]{On the partial sums of Walsh-Fourier series}
\address{G. Tephnadze, Department of Mathematics, Faculty of Exact and
Natural Sciences, Tbilisi State University, Chavchavadze str. 1, Tbilisi
0128, Georgia and Department of Engineering Sciences and Mathematics, Lule%
\aa\ University of Technology, SE-971 87 Lule\aa , Sweden.}
\email{giorgitephnadze@gmail.com}
\thanks{The research was supported by Shota Rustaveli National Science
Foundation grants no.13/06 and DO/24/5-100/14.}
\date{}
\maketitle

\begin{abstract}

In this paper we investigate some convergence and divergence of some specific subsequences of partial sums with respect to Walsh system on the martingale Hardy spaces. By using these results we obtain relationship of the ratio of convergence of the partial sums of the Walsh series with the modulus of continuity of martingale. These conditions are in a sense necessary and sufficient.

\end{abstract}

\date{}

\textbf{2010 Mathematics Subject Classification.} 42C10.

\textbf{Key words and phrases:} Walsh system, partial sums, martingale Hardy
space.

\section{INTRODUCTION}

It is well-known that (see e.g. \cite{1} and \cite{13}) Walsh system does
not form basis in the space $L_{1}.$ Moreover, there exists a function in
the dyadic Hardy space $H_{1},$ such that the partial sums of $f$ are not
bounded in $L_{1}$-norm, but partial sums $S_{n}$ of the Walsh-Fourier
series of a function $f\in L_{1}$ convergence in measure (see also \cite{ggt}
and \cite{gol}).

For Walsh-Fourier series Onneweer \cite{10} showed that if modulus of
continuity of $f\in L_{1}\left[ 0,1\right) $ satisfies the condition 
\begin{equation}
\omega _{1}\left( \delta ,f\right) =o\left( \frac{1}{\log \left( 1/\delta
\right) }\right)\text{ \ as \ \ }\delta \rightarrow 0,  \label{1aaa}
\end{equation}%
then its Walsh-Fourier series converges in $L_{1}$-norm. He also proved that
condition (\ref{1aaa}) can not be improved.

It is also known that subsequence of partial sums $S_{m_{k}}$ is bounded
from $L_{1}$ to $L_{1}$ if and only if $\left\{ m_{k}:k\geq 0\right\} $ have
uniformly bounded variations. In \cite[Ch. 1]{13} it was proved that if 
\textbf{\ }$f\in L_{1}\left( G\right) $ and $\left\{ m_{n}:n\geq 1\right\} $
be subsequence of positive numbers $\mathbb{N}$, such that 
\begin{equation}
\omega _{1}\left( \frac{1}{m_{n}},f\right) =o\left( \frac{1}{L_{S}\left(
m_{n}\right) }\right)\text{ \ as \ \ }n\rightarrow \infty ,  \label{2aaaa}
\end{equation}%
where the number $L_{S}\left( n\right) $ is $n$-th Lebesgue constant, then
subsequence $S_{m_{n}}f$ of partial sums of Walsh-Fourier series converge in 
$L_{1}$-norm. Goginava and Tkebuchava \cite{gt} proved that the condition (%
\ref{2aaaa}) can not be improved. Since (see \cite{luk} and e.g. \cite{sws}) 
\begin{equation}
\frac{V\left( n\right) }{8}\leq L_{S}\left( n\right) \leq V\left( n\right)
\label{22b}
\end{equation}%
the condition (\ref{2aaaa}) can be rewritten with the condition 
\begin{equation*}
\omega _{1}\left( \frac{1}{m_{n}},f\right) =o\left( \frac{1}{V\left(
m_{n}\right) }\right)\text{ \ as \ \ }n\rightarrow \infty .
\end{equation*}%
In \cite{tep1} it was proved that if $F\in H_{p}$ and 
\begin{equation}
\omega _{H_{p}}\left( \frac{1}{2^{n}},F\right) =o\left( \frac{1}{\left( n^{%
\left[ p\right] }2^{\left( 1/p-1\right) n}\right) }\right) \text{ as }%
n\rightarrow \infty ,  \label{3aaa}
\end{equation}%
where $0<p\leq 1$ and $\left[ p\right] $ denotes integer part of $p,$ then $%
S_{n}F\rightarrow F$ as $n\rightarrow \infty $ in $L_{p,\infty }$-norm.
Moreover, there was showed that condition (\ref{3aaa}) can not be improved.

Uniform and pointwise convergence and some approximation properties of
partial sums in $L_{1}$ norms was investigate by Goginava \cite{gog1} (see
also \cite{gt}, \cite{gog2}), Nagy \cite{na}, Avdispahić and Memić \cite{am}%
. Fine \cite{fi} obtained sufficient conditions for the uniform convergence
which are complete analogy with the Dini-Lipschits conditions. Guličev \cite%
{9} estimated the rate of uniform convergence of a Walsh-Fourier series by
using Lebesgue constants and modulus of continuity. These problems for
Vilenkin groups were considered by Blahota \cite{2}, Fridli \cite{4} and Gát 
\cite{5}.

The main aim of this paper is to find characterizations of boundedness of
the subsequence of partial sums of the Walsh series of $H_{p}$ martingales
in terms of measurable properties of a Dirichlet kernel corresponding to
partial summing. As a consequence we get the corollaries about the
convergence and divergence of some specific subsequences of partial sums.
For $p=1$ the simple numerical criterion for the index of partial sum in
terms of its dyadic expansion is given which governs the convergence (or the
ratio of divergence). Another type of results covered by the paper is the
relationship of the ratio of convergence of the partial sums of the Walsh
series with the modulus of continuity of martingale. The conditions given
below are in a sense necessary and sufficient.

\section{Preliminaries}

\bigskip Let $\mathbb{N}_{+}$ denote the set of the positive integers, $%
\mathbb{N}:=\mathbb{N}_{+}\cup \{0\}.$ Denote by $Z_{2}$ the discrete cyclic
group of order 2, that is $Z_{2}:=\{0,1\},$ where the group operation is the
modulo 2 addition and every subset is open. The Haar measure on $Z_{2}$ is
given so that the measure of a singleton is 1/2.

Define the group $G$ as the complete direct product of the group $Z_{2},$
with the product of the discrete topologies of $Z_{2}$`s. The elements of $G$
are represented by sequences $x:=(x_{0},x_{1},...,x_{j},...),$ where $%
x_{k}=0\vee 1.$

It is easy to give a base for the neighborhood of $x\in G$ 
\begin{equation*}
I_{0}\left( x\right) :=G,\text{ \ }I_{n}(x):=\{y\in
G:y_{0}=x_{0},...,y_{n-1}=x_{n-1}\}\text{ }(n\in \mathbb{N}).
\end{equation*}

Denote $I_{n}:=I_{n}\left( 0\right) ,$ $\overline{I_{n}}:=G$ $\backslash $ $%
I_{n}$ and $e_{n}:=\left( 0,...,0,x_{n}=1,0,...\right) \in G,$ for $n\in 
\mathbb{N}$. Then it is easy to show that

\begin{equation}
\overline{I_{M}}=\overset{M-1}{\underset{s=0}{\bigcup }}I_{s}\backslash
I_{s+1}.  \label{2}
\end{equation}

If $n\in \mathbb{N},$ then every $n$ can be uniquely expressed as $%
n=\sum_{k=0}^{\infty }n_{j}2^{j},$ where $n_{j}\in Z_{2}$ $~(j\in \mathbb{N}%
) $ and only a finite numbers of $n_{j}$ differ from zero.

Let 
\begin{equation*}
\left\langle n\right\rangle :=\min \{j\in \mathbb{N},n_{j}\neq 0\}\text{ \ \
and \ \ \ }\left\vert n\right\vert :=\max \{j\in \mathbb{N},n_{j}\neq 0\},
\end{equation*}%
that is $2^{\left\vert n\right\vert }\leq n\leq 2^{\left\vert n\right\vert
+1}.$ Set 
\begin{equation*}
d\left( n\right) =\left\vert n\right\vert -\left\langle n\right\rangle ,%
\text{ \ for \ all \ \ }n\in \mathbb{N}.
\end{equation*}

Define the variation of an $n\in \mathbb{N}$ with binary coefficients $%
\left( n_{k},\text{ }k\in \mathbb{N}\right) $ by

\begin{equation*}
V\left( n\right) =n_{0}+\overset{\infty }{\underset{k=1}{\sum }}\left|
n_{k}-n_{k-1}\right| .
\end{equation*}

The norms (or quasi-norm) of the spaces $L_{p}(G)$ and $L_{p,\infty }\left(
G\right) ,$ $\left( 0<p<\infty \right) $ are respectively defined by 
\begin{equation*}
\left\Vert f\right\Vert _{p}^{p}:=\int_{G}\left\vert f\right\vert ^{p}d\mu
,\ \ \ \ \left\Vert f\right\Vert _{L_{p,\infty }}^{p}:=\sup_{\lambda
>0}\lambda ^{p}\mu \left( f>\lambda \right) .
\end{equation*}

The $k$-th Rademacher function is defined by%
\begin{equation*}
r_{k}\left( x\right) :=\left( -1\right) ^{x_{k}}\text{\qquad }\left( \text{ }%
x\in G,\text{ }k\in \mathbb{N}\right) .
\end{equation*}

Now, define the Walsh system $w:=(w_{n}:n\in \mathbb{N})$ on $G$ as: 
\begin{equation*}
w_{n}(x):=\overset{\infty }{\underset{k=0}{\Pi }}r_{k}^{n_{k}}\left(
x\right) =r_{\left\vert n\right\vert }\left( x\right) \left( -1\right) ^{%
\underset{k=0}{\overset{\left\vert n\right\vert -1}{\sum }}n_{k}x_{k}}\text{%
\qquad }\left( n\in \mathbb{N}\right) .
\end{equation*}

The Walsh system is orthonormal and complete in $L_{2}\left( G\right) $ (see
e.g. \cite{sws}).

If $f\in L_{1}\left( G\right) $ we can establish the Fourier coefficients,
the partial sums of the Fourier series, the Dirichlet kernels with respect
to the Walsh system in the usual manner:$\qquad $%
\begin{equation*}
\widehat{f}\left( k\right) :=\int_{G}fw_{k}d\mu \,\,\,\,\left( k\in \mathbb{N%
}\right) ,
\end{equation*}%
\begin{equation*}
S_{n}f:=\sum_{k=0}^{n-1}\widehat{f}\left( k\right) w_{k},\ \
D_{n}:=\sum_{k=0}^{n-1}w_{k\text{ }}\,\,\,\left( n\in \mathbb{N}_{+}\right) .
\end{equation*}

Recall that%
\begin{equation}
D_{2^{n}}\left( x\right) =\left\{ 
\begin{array}{ll}
2^{n}, & \,\text{if\thinspace \thinspace \thinspace }x\in I_{n} \\ 
0, & \text{if}\,\,x\notin I_{n}%
\end{array}%
\right.  \label{1dn}
\end{equation}%
and

\begin{equation}
D_{n}=w_{n}\overset{\infty }{\underset{k=0}{\sum }}n_{k}r_{k}D_{2^{k}}=w_{n}%
\overset{\infty }{\underset{k=0}{\sum }}n_{k}\left(
D_{2^{k+1}}-D_{2^{k}}\right) ,\text{ for \ }n=\overset{\infty }{\underset{i=0%
}{\sum }}n_{i}2^{i}.  \label{2dn}
\end{equation}%
\ 

Let us denote $n$-th Lebesgue constant by%
\begin{equation*}
L_{S}\left( n\right) :=\left\Vert D_{n}\right\Vert _{1}.
\end{equation*}

The $\sigma $-algebra generated by the intervals $\left\{ I_{n}\left(
x\right) :x\in G\right\} $ will be denoted by $\zeta _{n}\left( n\in \mathbb{%
N}\right).$ Denote by $F=\left( F_{n},n\in \mathbb{N}\right) $ the
martingale with respect to $\digamma _{n}$ $\left( n\in \mathbb{N}\right) $
(for details see e.g. \cite{We1}).

The maximal function of a martingale $F$ is defined by

\begin{equation*}
F^{\ast }=\sup_{n\in \mathbb{N}}\left\vert F_{n}\right\vert .
\end{equation*}

In case $f\in L_{1}\left( G\right) ,$ the maximal functions are also be
given by

\begin{equation*}
f^{\ast }\left( x\right) =\sup\limits_{n\in \mathbb{N}}\frac{1}{\mu \left(
I_{n}\left( x\right) \right) }\left\vert \int\limits_{I_{n}\left( x\right)
}f\left( u\right) d\mu \left( u\right) \right\vert .
\end{equation*}

For $0<p<\infty $ the Hardy martingale spaces $H_{p}\left( G\right) $
consist all martingale for which

\begin{equation*}
\left\Vert F\right\Vert _{H_{p}}:=\left\Vert F^{\ast }\right\Vert
_{p}<\infty .
\end{equation*}

The best approximation of $f\in L_{p}(G)$ $(1\leq p\in \infty )$ is defined
as%
\begin{equation*}
E_{n}\left( f,L_{p}\right) =\inf_{\psi \in \emph{p}_{n}}\left\Vert f-\psi
\right\Vert _{p},
\end{equation*}%
where $\emph{p}_{n}$ is set of all Walsh polynomials of order less than $%
n\in \mathbb{N}$.

The integrated modulus of continuity of $f\in L_{p}$ is defined by

\begin{equation*}
\omega _{p}\left( \frac{1}{2^{n}},f\right) =\sup\limits_{h\in
I_{n}}\left\Vert f\left( \cdot +h\right) -f\left( \cdot \right) \right\Vert
_{p}.
\end{equation*}

The concept of modulus of continuity in $H_{p}\left( G\right) $ $\left(
0<p\leq 1\right) $ can be defined in following way 
\begin{equation*}
\omega _{H_{p}}\left( \frac{1}{2^{n}},F\right) :=\left\Vert
F-S_{2^{n}}F\right\Vert _{H_{p}}.
\end{equation*}

Watari \cite{wat} showed that there are strong connections between 
\begin{equation*}
\omega _{p}\left( \frac{1}{2^{n}},f\right) ,\text{ \ }E_{2^{n}}\left(
f,L_{p}\right) \text{ \ \ and\ \ \ }\left\Vert f-S_{2^{n}}f\right\Vert
_{p},\ \ p\geq 1,\text{ \ }n\in \mathbb{N}.
\end{equation*}

In particular,%
\begin{equation}
\frac{1}{2}\omega _{p}\left( \frac{1}{2^{n}},f\right) \leq \left\Vert
f-S_{2^{n}}f\right\Vert _{p}\leq \omega _{p}\left( \frac{1}{2^{n}},f\right) ,
\label{eqvi}
\end{equation}%
and%
\begin{equation*}
\frac{1}{2}\left\Vert f-S_{2^{n}}f\right\Vert _{p}\leq E_{2^{n}}\left(
f,L_{p}\right) \leq \left\Vert f-S_{2^{n}}f\right\Vert _{p}.
\end{equation*}

A bounded measurable function $a$ is called p-atom, if there exist a dyadic
interval $I,$ such that \qquad 
\begin{equation*}
\int_{I}ad\mu =0,\text{ \ \ }\left\Vert a\right\Vert _{\infty }\leq \mu
\left( I\right) ^{-1/p},\text{ \ \ supp}\left( a\right) \subset I.
\end{equation*}

The dyadic Hardy martingale spaces $H_{p}$ for $0<p\leq 1$ have an atomic
characterization. Namely, the following theorem is true (see \cite{S} and 
\cite{We3}):

\textbf{Theorem W: }A martingale $F=\left( F_{n},n\in \mathbb{N}\right) $ is
in $H_{p}\left( 0<p\leq 1\right) $ if and only if there exists a sequence $%
\left( a_{k},\text{ }k\in \mathbb{N}\right) $ of p-atoms and a sequence $%
\left( \mu _{k},k\in \mathbb{N}\right) $ of a real numbers such that for
every $n\in \mathbb{N}$

\begin{equation}
\qquad \sum_{k=0}^{\infty }\mu _{k}S_{2^{n}}a_{k}=F_{n},\text{ \ \ \ }%
\sum_{k=0}^{\infty }\left\vert \mu _{k}\right\vert ^{p}<\infty ,  \label{2A}
\end{equation}

Moreover, 
\begin{equation*}
\left\Vert F\right\Vert _{H_{p}}\backsim \inf \left( \sum_{k=0}^{\infty
}\left\vert \mu _{k}\right\vert ^{p}\right) ^{1/p},
\end{equation*}
where the infimum is taken over all decomposition of $F$ of the form (\ref%
{2A}).

It is easy to check that for every martingale $F=\left( F_{n},n\in \mathbb{N}%
\right) $ and every $k\in \mathbb{N}$ the limit

\begin{equation}
\widehat{F}\left( k\right) :=\lim_{n\rightarrow \infty }\int_{G}F_{n}\left(
x\right) w_{k}\left( x\right) d\mu \left( x\right)  \label{3a}
\end{equation}%
exists and it is called the $k$-th Walsh-Fourier coefficients of $F.$

If $F:=$ $\left( E_{n}f:n\in \mathbb{N}\right) $ is regular martingale,
generated by $f\in L_{1}\left( G\right) ,$ then \ $\widehat{F}\left(
k\right) =\widehat{f}\left( k\right) ,$ $k\in \mathbb{N}.$

For the martingale 
\begin{equation*}
F=\sum_{n=0}^{\infty }\left( F_{n}-F_{n-1}\right)
\end{equation*}%
\ the conjugate transforms are defined as 
\begin{equation*}
\widetilde{F^{\left( t\right) }}=\sum_{n=0}^{\infty }r_{n}\left( t\right)
\left( F_{n}-F_{n-1}\right) ,
\end{equation*}%
where $t\in G$ is fixed. Note that $\widetilde{F^{\left( 0\right) }}=F.$ As
is well known (see e.g. \cite{We1}) 
\begin{equation}
\left\Vert \widetilde{F^{\left( t\right) }}\right\Vert _{H_{p}}=\left\Vert
F\right\Vert _{H_{p}},\text{ \ }\left\Vert F\right\Vert _{H_{p}}^{p}\sim
\int_{G}\left\Vert \widetilde{F^{\left( t\right) }}\right\Vert _{p}^{p}dt,%
\text{ \ \ }\widetilde{S_{n}F^{\left( t\right) }}=S_{n}\widetilde{F^{\left(
t\right) }}.  \label{5.1}
\end{equation}

\section{Formulation of Main Results}

\begin{theorem}
a) Let $0<p<1$ and $F\in H_{p}$. Then there exists an \textit{absolute
constant }$c_{p},$ depending only on $p,$ such that%
\begin{equation*}
\text{ }\left\Vert S_{n}F\right\Vert _{H_{p}}\leq c_{p}2^{d\left( n\right)
\left( 1/p-1\right) }\left\Vert F\right\Vert _{H_{p}}.
\end{equation*}

\textit{b) Let} $0<p<1,$ $\left\{ m_{k}:\text{ }k\geq 0\right\} $ \textit{be
any increasing sequence of positive integers} $\mathbb{N}_{+}$ \textit{such
that }%
\begin{equation}
\sup_{k\in \mathbb{N}}d\left( m_{k}\right) =\infty  \label{dnk}
\end{equation}%
\ \textit{and } $\Phi :\mathbb{N}_{+}\rightarrow \lbrack 1,\infty )$ \textit{%
be any nondecreasing function, satisfying condition} 
\begin{equation}
\overline{\underset{k\rightarrow \infty }{\lim }}\frac{2^{d\left(
m_{k}\right) \left( 1/p-1\right) }}{\Phi \left( m_{k}\right) }=\infty .
\label{1010}
\end{equation}%
\textit{Then there exists a martingale }$F\in H_{p},$ \textit{such that} 
\begin{equation*}
\underset{k\in \mathbb{N}}{\sup }\left\Vert \frac{S_{m_{k}}F}{\Phi \left(
m_{k}\right) }\right\Vert _{L_{p,\infty }}=\infty .
\end{equation*}
\end{theorem}

\begin{corollary}
a) Let $0<p<1$ and $F\in H_{p}$. Then there exists an \textit{absolute
constant }$c_{p},$ depending only on $p,$ such that%
\begin{equation*}
\text{ }\left\Vert S_{n}F\right\Vert _{H_{p}}\leq c_{p}\left( n\mu \left\{ 
\text{supp}\left( D_{n}\right) \right\} \right) ^{1/p-1}\left\Vert
F\right\Vert _{H_{p}}.
\end{equation*}

\textit{b) Let} $0<p<1,$ $\left\{ m_{k}:\text{ }k\geq 0\right\} $ \textit{be
any increasing sequence of positive integers} $\mathbb{N}_{+}$ \textit{such
that }%
\begin{equation}
\sup_{k\in \mathbb{N}}m_{k}\mu \left\{ \text{supp}\left( D_{m_{k}}\right)
\right\} =\infty  \label{suppdnk}
\end{equation}%
\ \textit{and } $\Phi :\mathbb{N}_{+}\rightarrow \lbrack 1,\infty )$ \textit{%
be any nondecreasing function, satisfying condition} 
\begin{equation}
\overline{\underset{k\rightarrow \infty }{\lim }}\frac{\left( m_{k}\mu
\left\{ \text{supp}\left( D_{m_{k}}\right) \right\} \right) ^{1/p-1}}{\Phi
\left( m_{k}\right) }=\infty .  \label{12e}
\end{equation}%
\textit{Then there exists a martingale }$F\in H_{p},$ \textit{such that} 
\begin{equation*}
\underset{k\in \mathbb{N}}{\sup }\left\Vert \frac{S_{m_{k}}F}{\Phi \left(
m_{k}\right) }\right\Vert _{L_{p,\infty }}=\infty .
\end{equation*}
\end{corollary}

\begin{corollary}
Let $n\in \mathbb{N}$ and $0<p<1.$ Then there exists a martingale $F\in
H_{p},$ such that 
\begin{equation}
\underset{n\in \mathbb{N}}{\sup }\left\Vert S_{2^{n}+1}F\right\Vert
_{L_{p,\infty }}=\infty .  \label{sn2n1}
\end{equation}
\end{corollary}

\begin{corollary}
Let $n\in \mathbb{N}$ and $0<p\leq 1$ and $F\in H_{p}$. Then 
\begin{equation}
\left\Vert S_{2^{n}+2^{n-1}}F\right\Vert _{H_{p}}\leq c_{p}\left\Vert
F\right\Vert _{H_{p}}.  \label{sn2n2}
\end{equation}
\end{corollary}

\begin{theorem}
a) Let $n\in \mathbb{N}_{+}$ and $F\in H_{1}.$ Then there exists an absolute
constant $c,$ such that 
\begin{equation*}
\left\Vert S_{n}F\right\Vert _{H_{1}}\leq cV\left( n\right) \left\Vert
F\right\Vert _{H_{1}}.
\end{equation*}

\textit{b) Let} $\left\{ m_{k}:\text{ }k\geq 0\right\} $ \textit{be any
increasing sequence of positive integers} $\mathbb{N}_{+},$ \textit{such that%
} 
\begin{equation}
\sup_{k\in \mathbb{N}}V\left( m_{k}\right) =\infty  \label{vnk}
\end{equation}%
\textit{and} $\Phi :\mathbb{N}_{+}\rightarrow \lbrack 1,\infty )$ \textit{be
any nondecreasing function, satisfying condition} 
\begin{equation}
\overline{\underset{k\rightarrow \infty }{\lim }}\frac{V\left( m_{k}\right) 
}{\Phi \left( m_{k}\right) }=\infty .  \label{17aa}
\end{equation}

\textit{Then there exist a martingale} $F\in H_{1},$ \textit{such that} 
\begin{equation*}
\underset{{k\in \mathbb{N}}}{\sup }\left\Vert \frac{S_{m_{k}}F}{\Phi \left(
m_{k}\right) }\right\Vert _{1}=\infty .
\end{equation*}
\end{theorem}

\begin{theorem}
Let $2^{k}<n\leq 2^{k+1}.$ Then there exist absolute constant $c_{p},$
depending only on $p,$ such that 
\begin{equation*}
\left\Vert S_{n}F-F\right\Vert _{H_{p}}\leq c_{p}2^{d\left( n\right) \left(
1/p-1\right) }\omega _{H_{p}}\left( \frac{1}{2^{k}},F\right) ,\text{ \ \ \ }%
\left( 0<p<1\right)
\end{equation*}%
and%
\begin{equation}
\left\Vert S_{n}F-F\right\Vert _{H_{1}}\leq c_{1}V\left( n\right) \omega
_{H_{1}}\left( \frac{1}{2^{k}},F\right) .  \label{sn2}
\end{equation}
\end{theorem}

\begin{theorem}
a) \textit{Let }$0<p<1,$ $F\in H_{p}$ and $\{m_{k}:k\geq 0\}$\textit{\ be a
sequence of nonnegative integers},\textit{\ such that} 
\begin{equation}
\omega _{H_{p}}\left( \frac{1}{2^{\left\vert m_{k}\right\vert }},F\right)
=o\left( \frac{1}{2^{d\left( m_{k}\right) \left( 1/p-1\right) }}\right) 
\text{ as \ }k\rightarrow \infty .  \label{18a}
\end{equation}%
Then 
\begin{equation}
\left\Vert S_{m_{k}}F-F\right\Vert _{H_{p}}\rightarrow 0\text{ \ as \ }%
k\rightarrow \infty .  \label{con1}
\end{equation}

\textit{b) Let} $\{m_{k}:k\geq 0\}$ \textit{be any increasing sequence of
positive integers} $\mathbb{N}_{+},$ \textit{satisfying condition (\ref{dnk}%
).} \textit{Then there exists a martingale} $F\in H_{p}$ \textit{and
subsequence} $\{\alpha _{k}:k\geq 0\}\subset \{m_{k}:k\geq 0\},$ \textit{for
which} 
\begin{equation*}
\omega _{H_{p}}\left( \frac{1}{2^{\left\vert \alpha _{k}\right\vert }}%
,F\right) =O\left( \frac{1}{2^{d\left( \alpha _{k}\right) \left(
1/p-1\right) }}\right) \text{ as \ }k\rightarrow \infty \text{\ }
\end{equation*}%
\textit{and} 
\begin{equation}
\limsup\limits_{k\rightarrow \infty }\left\Vert S_{\alpha
_{k}}F-F\right\Vert _{L_{p,\infty }}>c_{p}>0\,\,\,\text{as\thinspace
\thinspace \thinspace }k\rightarrow \infty ,  \label{con11}
\end{equation}%
where $c_{p}$ is an absolute constant depending only on $p.$
\end{theorem}

\begin{corollary}
a) \textit{Let }$0<p<1,$ $F\in H_{p}$ and $\{m_{k}:k\geq 0\}$\textit{\ be a
sequence of nonnegative integers},\textit{\ such that} 
\begin{equation}
\omega _{H_{p}}\left( \frac{1}{2^{\left\vert m_{k}\right\vert }},F\right)
=o\left( \frac{1}{\left( m_{k}\mu \left( \text{supp}D_{m_{k}}\right) \right)
^{1/p-1}}\right) \text{ as \ }k\rightarrow \infty .  \label{cond2}
\end{equation}%
Then (\ref{con1}) is satisfied.

\textit{b) Let} $\{m_{k}:k\geq 0\}$ \textit{be any increasing sequence of
positive integers} $\mathbb{N}_{+},$ \textit{satisfying condition (\ref%
{suppdnk}).} \textit{Then there exists a martingale} $F\in H_{p}$ \textit{%
and subsequence} $\{\alpha _{k}:k\geq 0\}\subset \{m_{k}:k\geq 0\},$ \textit{%
for which} 
\begin{equation*}
\omega _{H_{p}}\left( \frac{1}{2^{\left\vert \alpha _{k}\right\vert }}%
,F\right) =O\left( \frac{1}{\left( \alpha _{k}\mu \left( \text{supp}%
D_{\alpha _{k}}\right) \right) ^{1/p-1}}\right) \text{\ \ as \ }k\rightarrow
\infty .
\end{equation*}%
\textit{and }(\ref{con11}) is satisfied.
\end{corollary}

\begin{theorem}
a) \textit{Let }$F\in H_{1}$ and $\{m_{k}:k\geq 0\}$\textit{\ be a sequence
of nonnegative integers},\textit{\ such that} 
\begin{equation}
\omega _{H_{1}}\left( \frac{1}{2^{\left\vert m_{k}\right\vert }},F\right)
=o\left( \frac{1}{V\left( m_{k}\right) }\right) \text{ as \ }k\rightarrow
\infty .  \label{cond1}
\end{equation}%
Then 
\begin{equation}
\left\Vert S_{m_{k}}F-F\right\Vert _{H_{1}}\rightarrow 0\text{ \ as \ }%
k\rightarrow \infty .  \label{cond1a}
\end{equation}

b) \textit{Let }$\{m_{k}:k\geq 0\}$\textit{\ be any increasing sequence of
positive integers} $\mathbb{N}_{+}$,\textit{\ satisfying condition (\ref{vnk}%
).} \textit{Then there exists a martingale }$F\in H_{1}$\textit{\ and }$%
\{\alpha _{k}:k\geq 0\}\subset \{m_{k}:k\geq 0\}$\textit{\ for which} 
\begin{equation*}
\omega _{H_{1}}\left( \frac{1}{2^{\left\vert \alpha _{k}\right\vert }}%
,F\right) =O\left( \frac{1}{V\left( \alpha _{k}\right) }\right) \text{ \ as
\ }k\rightarrow \infty .
\end{equation*}%
\textit{and} 
\begin{equation}
\limsup\limits_{k\rightarrow \infty }\left\Vert S_{\alpha
_{k}}F-F\right\Vert _{1}>c>0\,\,\,\text{as\thinspace \thinspace \thinspace }%
k\rightarrow \infty ,  \label{cond10}
\end{equation}%
where $c$ is an absolute constant.
\end{theorem}

\section{Proof of the Theorems}

\begin{proof}[\textbf{Proof of Theorem 1.}]
\textbf{\ \ }Suppose that 
\begin{equation}
\left\Vert 2^{\left( 1-1/p\right) d\left( n\right) }S_{n}F\right\Vert
_{p}\leq c_{p}\left\Vert F\right\Vert _{H_{p}}.  \label{11.1}
\end{equation}

By combining (\ref{5.1}) and (\ref{11.1}) we get that 
\begin{equation}
\left\Vert 2^{\left( 1-1/p\right) d\left( n\right) }S_{n}F\right\Vert
_{H_{p}}\leq c_{p}\int_{G}\left\Vert 2^{\left( 1-1/p\right) d\left( n\right)
}\widetilde{S_{n}F^{\left( t\right) }}\right\Vert _{p}d\mu \left( t\right)
\label{11.2}
\end{equation}%
\begin{equation*}
=c_{p}\int_{G}\left\Vert 2^{\left( 1-1/p\right) d\left( n\right) }S_{n}%
\widetilde{F^{\left( t\right) }}\right\Vert _{p}d\mu \left( t\right) \leq
c_{p}\int_{G}\left\Vert \widetilde{F^{\left( t\right) }}\right\Vert
_{H_{p}}d\mu \left( t\right) \leq c_{p}\left\Vert F\right\Vert _{H_{p}}.
\end{equation*}

By using Theorem W, (\ref{11.2}) the proof of Theorem 1 will be complete, if
we show that%
\begin{equation}
\int\limits_{G}\left\vert 2^{\left( 1-1/p\right) d\left( n\right)
}S_{n}a\right\vert ^{p}d\mu \leq c_{p}<\infty ,  \label{25a}
\end{equation}%
for every p-atom $a,$ with support$\ I$ and $\mu \left( I\right) =2^{-N}$.

We may assume that this arbitrary $p$-atom $a$ has support $I=I_{M}.$ It is
easy to see that $S_{n}\left( a\right) =0,$ when $2^{M}$ $\geq n$.
Therefore, we can suppose that $2^{M}<n$. Since $\left\Vert a\right\Vert
_{\infty }\leq 2^{M/p}$ we can write that 
\begin{equation}
\left\vert 2^{\left( 1-1/p\right) d\left( n\right) }S_{n}a\left( x\right)
\right\vert \leq 2^{\left( 1-1/p\right) d\left( n\right) }\left\Vert
a\right\Vert _{\infty }\int_{I_{M}}\left\vert D_{n}\left( x+t\right)
\right\vert d\mu \left( t\right)  \label{11a}
\end{equation}%
\begin{equation*}
\leq 2^{M/p}2^{\left( 1-1/p\right) d\left( n\right) }\int_{I_{M}}\left\vert
D_{n}\left( x+t\right) \right\vert d\mu \left( t\right) .
\end{equation*}

Let $x\in I_{M}$. Since $V\left( n\right) \leq d\left( n\right) ,$ by
applying (\ref{22b}) we get that 
\begin{equation*}
\left\vert 2^{\left( 1-1/p\right) d\left( n\right) }S_{n}a\right\vert \leq
2^{M/p}2^{\left( 1-1/p\right) d\left( n\right) }V\left( n\right) \leq
2^{M/p}d\left( n\right) 2^{\left( 1-1/p\right) d\left( n\right) }
\end{equation*}%
and 
\begin{equation}
\int_{I_{M}}\left\vert 2^{\left( 1-1/p\right) d\left( n\right)
}S_{n}a\right\vert ^{p}d\mu \leq d\left( n\right) 2^{\left( 1-1/p\right)
d\left( n\right) }<c_{p}<\infty .  \label{11b}
\end{equation}

Let $t\in I_{M}$ and $x\in I_{s}\backslash I_{s+1},\,0\leq s\leq
M-1<\left\langle n\right\rangle $ or $0\leq s<\left\langle n\right\rangle
\leq M-1.$ Then $x+t$ $\in I_{s}\backslash I_{s+1}$. By using (\ref{2dn}) we
get that $D_{n}\left( x+t\right) =0$ and 
\begin{equation*}
\left\vert 2^{\left( 1-1/p\right) d\left( n\right) }S_{n}a\left( x\right)
\right\vert =0.
\end{equation*}

Let $x\in I_{s}\backslash I_{s+1},$ $\,\left\langle n\right\rangle \leq
s\leq M-1.$ Then $x+t\in I_{s}\backslash I_{s+1},$ for $t\in I_{M}$. By
using (\ref{2dn}) we can write that 
\begin{equation*}
\left\vert D_{n}\left( x+t\right) \right\vert \leq
\sum_{j=0}^{s}n_{j}2^{j}\leq c2^{s}.
\end{equation*}%
If we apply (\ref{11a}) we get that%
\begin{equation}
\left\vert 2^{\left( 1-1/p\right) d\left( n\right) }S_{n}a\left( x\right)
\right\vert \leq 2^{\left( 1-1/p\right) d\left( n\right) }2^{M/p}\frac{2^{s}%
}{2^{M}}=2^{\left\langle n\right\rangle \left( 1/p-1\right) }2^{s}.
\label{12}
\end{equation}

By combining (\ref{2}) and (\ref{12}) we have 
\begin{equation*}
\int_{\overline{I_{M}}}\left\vert 2^{\left( 1-1/p\right) d\left( n\right)
}S_{n}a\left( x\right) \right\vert ^{p}d\mu \left( x\right)
\end{equation*}%
\begin{equation*}
=\overset{M-1}{\underset{s=\left\langle n\right\rangle }{\sum }}%
\int_{I_{s}\backslash I_{s+1}}\left\vert 2^{\left\langle n\right\rangle
\left( 1/p-1\right) }2^{s}\right\vert ^{p}d\mu \left( x\right) \leq c\overset%
{M-1}{\underset{s=\left\langle n\right\rangle }{\sum }}\frac{2^{\left\langle
n\right\rangle \left( 1-p\right) }}{2^{s\left( 1-p\right) }}\leq
c_{p}<\infty .
\end{equation*}

Let prove second part of Theorem 1. Under condition (\ref{1010}), there
exists sequence $\left\{ \alpha _{k}:\text{ }k\geq 0\right\} \subset \left\{
m_{k}:\text{ }k\geq 0\right\} ,$ such that 
\begin{equation}
\sum_{\eta =0}^{\infty }\frac{\Phi ^{p/2}\left( \alpha _{\eta }\right) }{%
2^{d\left( \alpha _{\eta }\right) \left( 1-p\right) /2}}<\infty ,
\label{12f}
\end{equation}

Let \qquad 
\begin{equation*}
F_{n}=\sum_{\left\{ k;\text{ }\left\vert \alpha _{k}\right\vert <n\right\}
}\lambda _{k}a_{k},\text{\ \ }
\end{equation*}%
where%
\begin{equation}
\lambda _{k}=\frac{\Phi ^{1/2}\left( \alpha _{k}\right) }{2^{d\left( \alpha
_{k}\right) \left( 1/p-1\right) /2}},\text{ \ \ }a_{k}=2^{\left\vert \alpha
_{k}\right\vert \left( 1/p-1\right) }\left( D_{2^{\left\vert \alpha
_{k}\right\vert +1}}-D_{2^{\left\vert \alpha _{k}\right\vert }}\right) ,
\label{100}
\end{equation}

By combining Theorem W and (\ref{12f}) we conclude that $F=\left( F_{n},n\in 
\mathbb{N}\right) \in H_{p}.$

By simple calculation we get\ that\thinspace 
\begin{equation}
\widehat{F}(j)  \label{6}
\end{equation}%
\begin{equation*}
=\left\{ 
\begin{array}{c}
\Phi ^{1/2}\left( \alpha _{k}\right) 2^{\left( \left\vert \alpha
_{k}\right\vert +\left\langle \alpha _{k}\right\rangle \right) \left(
1/p-1\right) /2},\text{ \ if \thinspace \thinspace }j\in \left\{
2^{_{\left\vert \alpha _{k}\right\vert }},...,2^{_{\left\vert \alpha
_{k}\right\vert +1}}-1\right\} ,\text{ }k=0,1,... \\ 
0,\text{ \ if \thinspace \thinspace \thinspace }j\notin
\bigcup\limits_{k=0}^{\infty }\left\{ 2^{_{\left\vert \alpha _{k}\right\vert
}},...,2^{_{\left\vert \alpha _{k}\right\vert +1}}-1\right\} .\text{ }%
\end{array}%
\right.
\end{equation*}

Since 
\begin{equation}
D_{j+2^{n}}=D_{2^{n}}+w_{2^{n}}D_{j},\text{ when }\,\,j\leq 2^{n},
\label{6aa}
\end{equation}%
by applying (\ref{6}) we have%
\begin{equation}
\frac{S_{\alpha _{k}}F}{\Phi \left( \alpha _{k}\right) }=\frac{1}{\Phi
\left( \alpha _{k}\right) }\sum_{\eta =0}^{k-1}\sum_{v=2^{\left\vert \alpha
_{\eta }\right\vert }}^{2^{\left\vert \alpha _{\eta }\right\vert +1}-1}%
\widehat{F}(v)w_{v}+\frac{1}{\Phi \left( \alpha _{k}\right) }%
\sum_{v=2^{\left\vert \alpha _{k}\right\vert }}^{\alpha _{k}-1}\widehat{F}%
(v)w_{v}  \label{6aaa}
\end{equation}%
\begin{equation*}
=\frac{1}{\Phi \left( \alpha _{k}\right) }\sum_{\eta
=0}^{k-1}\sum_{v=2^{\left\vert \alpha _{\eta }\right\vert }}^{2^{\left\vert
\alpha _{\eta }\right\vert +1}-1}\Phi ^{1/2}\left( \alpha _{\eta }\right)
2^{\left( \left\vert \alpha _{\eta }\right\vert +\left\langle \alpha _{\eta
}\right\rangle \right) \left( 1/p-1\right) /2}w_{v}
\end{equation*}%
\begin{equation*}
+\frac{1}{\Phi \left( \alpha _{k}\right) }\sum_{v=2^{\left\vert \alpha
_{k}\right\vert }}^{\alpha _{k}-1}\Phi ^{1/2}\left( \alpha _{k}\right)
2^{\left( \left\vert \alpha _{k}\right\vert +\left\langle \alpha
_{k}\right\rangle \right) \left( 1/p-1\right) /2}w_{v}
\end{equation*}%
\begin{equation*}
=\frac{1}{\Phi \left( \alpha _{k}\right) }\sum_{\eta =0}^{k-1}\Phi
^{1/2}\left( \alpha _{\eta }\right) 2^{\left( \left\vert \alpha _{\eta
}\right\vert +\left\langle \alpha _{\eta }\right\rangle \right) \left(
1/p-1\right) /2}\left( D_{2^{\left\vert \alpha _{\eta }\right\vert
+1}}-D_{2^{\left\vert \alpha _{\eta }\right\vert +1}}\right)
\end{equation*}%
\begin{equation*}
+\frac{2^{\left( \left\vert \alpha _{k}\right\vert +\left\langle \alpha
_{k}\right\rangle \right) \left( 1/p-1\right) /2}w_{2^{\left\vert \alpha
_{k}\right\vert }}D_{\alpha _{k}-2^{\left\vert \alpha _{k}\right\vert }}}{%
\Phi ^{1/2}\left( \alpha _{k}\right) }:=I+II.
\end{equation*}

By using (\ref{12f}) for $I$ we can write that 
\begin{equation}
\left\Vert I\right\Vert _{L_{p,\infty }}^{p}  \label{8a}
\end{equation}%
\begin{equation*}
\leq \frac{1}{\Phi ^{p}\left( \alpha _{k}\right) }\sum_{\eta =0}^{k-1}\frac{%
\Phi ^{p/2}\left( \alpha _{\eta }\right) }{2^{d\left( \alpha _{\eta }\right)
\left( 1-p\right) /2}}\left\Vert 2^{\left\vert \alpha _{\eta }\right\vert
\left( 1/p-1\right) }\left( D_{2^{\left\vert \alpha _{\eta }\right\vert
+1}}-D_{2^{\left\vert \alpha _{\eta }\right\vert +1}}\right) \right\Vert
_{L_{p,\infty }}^{p}
\end{equation*}%
\begin{equation*}
\leq \frac{1}{\Phi ^{p}\left( \alpha _{k}\right) }\sum_{\eta =0}^{\infty }%
\frac{\Phi ^{p/2}\left( \alpha _{\eta }\right) }{2^{d\left( \alpha _{\eta
}\right) \left( 1-p\right) /2}}\leq c<\infty .
\end{equation*}

Let $x\in I_{\left\langle \alpha _{k}\right\rangle }\backslash
I_{\left\langle \alpha _{k}\right\rangle +1}.$\textbf{\ }Under condition (%
\ref{2dn}) we can show that $\left\vert \alpha _{k}\right\vert \neq
\left\langle \alpha _{k}\right\rangle .$ It follows that $\left\langle
\alpha _{k}-2^{\left\vert \alpha _{k}\right\vert }\right\rangle
=\left\langle \alpha _{k}\right\rangle .\ $By combining (\ref{1dn}) and (\ref%
{2dn}) we have 
\begin{equation}
\left\vert D_{\alpha _{k}-2^{\left\vert \alpha _{k}\right\vert }}\right\vert
\label{77}
\end{equation}%
\begin{equation*}
=\left\vert \left( D_{2^{\left\langle \alpha _{k}\right\rangle
+1}}-D_{2^{\left\langle \alpha _{k}\right\rangle }}\right) +\overset{%
\left\vert \alpha _{k}\right\vert -1}{\underset{j=\left\langle \alpha
_{k}\right\rangle +1}{\sum }}\left( \alpha _{k}\right) _{j}\left(
D_{2^{i+1}}-D_{2^{i}}\right) \right\vert =\left\vert -D_{2^{\left\langle
\alpha _{k}\right\rangle }}\right\vert =2^{\left\langle \alpha
_{k}\right\rangle }
\end{equation*}%
and 
\begin{equation}
\left\vert II\right\vert =\frac{2^{\left( \left\vert \alpha _{k}\right\vert
+\left\langle \alpha _{k}\right\rangle \right) \left( 1/p-1\right) /2}}{\Phi
^{1/2}\left( \alpha _{k}\right) }\left\vert D_{\alpha _{k}-2^{\left\vert
\alpha _{k}\right\vert }}\left( x\right) \right\vert  \label{12aa}
\end{equation}%
\begin{equation*}
=\frac{2^{\left\vert \alpha _{k}\right\vert \left( 1/p-1\right)
/2}2^{\left\langle \alpha _{k}\right\rangle \left( 1/p+1\right) /2}}{\Phi
^{1/2}\left( \alpha _{k}\right) }.
\end{equation*}

By using (\ref{8a}) we see that 
\begin{equation*}
\left\Vert \frac{S_{\alpha _{k}}F}{\Phi \left( \alpha _{k}\right) }%
\right\Vert _{L_{p,\infty }}^{p}\geq \left\Vert II\right\Vert _{L_{p,\infty
}}^{p}-\left\Vert I\right\Vert _{L_{p,\infty }}^{p}
\end{equation*}%
\begin{equation*}
\geq \frac{2^{\left( \left\vert \alpha _{k}\right\vert \right) \left(
1/p-1\right) /2}2^{\left\langle \alpha _{k}\right\rangle \left( 1/p+1\right)
/2}}{\Phi ^{1/2}\left( \alpha _{k}\right) }\mu \left\{ x\in G:\text{ }%
\left\vert II\right\vert \geq \frac{2^{\left( \left\vert \alpha
_{k}\right\vert \right) \left( 1/p-1\right) /2}2^{\left\langle \alpha
_{k}\right\rangle \left( 1/p+1\right) /2}}{\Phi ^{1/2}\left( \alpha
_{k}\right) }\right\} ^{1/p}
\end{equation*}%
\begin{equation*}
\geq \frac{2^{\left( \left\vert \alpha _{k}\right\vert \right) \left(
1/p-1\right) /2}2^{\left\langle \alpha _{k}\right\rangle \left( 1/p+1\right)
/2}}{\Phi ^{1/2}\left( \alpha _{k}\right) }\left( \mu \left\{
I_{\left\langle \alpha _{k}\right\rangle }\backslash I_{\left\langle \alpha
_{k}\right\rangle +1}\right\} \right) ^{1/p}
\end{equation*}%
\begin{equation*}
\geq c\frac{2^{d\left( \alpha _{k}\right) \left( 1/p-1\right) /2}}{\Phi
^{1/2}\left( \alpha _{k}\right) }\rightarrow \infty ,\text{ \ as \ \ }%
k\rightarrow \infty .
\end{equation*}

Theorem 1 is proved.
\end{proof}

\begin{proof}[\textbf{Proof of Corollaries 1-3. }]
\bigskip By combining (\ref{1dn}) and (\ref{2dn}) we obtain 
\begin{equation*}
I_{\left\langle n\right\rangle }\backslash I_{\left\langle n\right\rangle
+1}\subset \text{supp}\left\{ D_{n}\right\} \subset I_{\left\langle
n\right\rangle }\text{ \ and }2^{-\left\langle n\right\rangle -1}\leq \mu
\left\{ \text{supp}\left( D_{n}\right) \right\} \leq 2^{-\left\langle
n\right\rangle }.\text{\ }
\end{equation*}%
It follows that 
\begin{equation*}
\frac{2^{d\left( n\right) \left( 1/p-1\right) }}{4}\leq \left( n\mu \left\{ 
\text{supp}\left( D_{n}\right) \right\} \right) ^{1/p-1}\leq 2^{d\left(
n\right) \left( 1/p-1\right) }.
\end{equation*}%
Corollary 1 is proved.

To prove Corollary 2 we only have to calculate that 
\begin{equation}
\left\vert 2^{n}+1\right\vert =n,\left\langle 2^{n}+1\right\rangle =0\text{
\ and \ }\rho \left( 2^{n}+1\right) =n.  \label{cor1}
\end{equation}%
By using the second part of Theorem 1 we see that there exists an martingale 
$F=\left( F_{n},n\in \mathbb{N}\right) \in H_{p}$ $\left( 0<p<1\right) ,$
such that (\ref{sn2n1}) is satisfied.

Let prove Corollary 3. \ Analogously to (\ref{cor1}) we can write that 
\begin{equation*}
\left\vert 2^{n}+2^{n-1}\right\vert =n,\left\langle
2^{n}+2^{n-1}\right\rangle =n-1\text{ \ and \ }\rho \left(
2^{n}+2^{n-1}\right) =1.
\end{equation*}%
By using the first part of Theorem 1 we immediately get inequality (\ref%
{sn2n2}), for all $0<p\leq 1$.

Corollaries 1-3 are proved.
\end{proof}

\begin{proof}[\textbf{Proof of Theorem 2. }]
By using (\ref{22b}) we have 
\begin{equation}
\left\Vert \frac{S_{n}F}{V\left( n\right) }\right\Vert _{1}\leq \left\Vert
F\right\Vert _{1}\leq \left\Vert F\right\Vert _{H_{1}}.  \label{12k}
\end{equation}

By combining (\ref{5.1}) and (\ref{12k}), after similar steps of (\ref{11.2}%
) we see that 
\begin{equation}
\left\Vert \frac{S_{n}F}{V\left( n\right) }\right\Vert _{H_{1}}\sim
\int_{G}\left\Vert \frac{\widetilde{S_{n}F^{\left( t\right) }}}{V\left(
n\right) }\right\Vert _{1}d\mu \left( t\right) \leq \left\Vert F\right\Vert
_{H_{1}}.  \label{50}
\end{equation}

Now, we prove second part of Theorem 2.\textbf{\ }Let $\left\{ m_{k}:\text{ }%
k\geq 0\right\} $ be subsequence of positive integers and function $\ \Phi :%
\mathbb{N}_{+}\rightarrow \lbrack 1,\infty )$ satisfies conditions of
Theorem 2. By using (\ref{17aa}) there exists an increasing sequence $%
\left\{ \alpha _{k}:\text{ }k\geq 0\right\} \subset \left\{ m_{k}:\text{ }%
k\geq 0\right\} $ of the positive integers such that 
\begin{equation}
\sum_{k=1}^{\infty }\frac{\Phi ^{1/2}\left( \alpha _{k}\right) }{%
V^{1/2}\left( \alpha _{k}\right) }\leq \beta <\infty .  \label{2aaa}
\end{equation}

Let \qquad 
\begin{equation*}
F_{n}:=\sum_{\left\{ k:\text{ }\left\vert \alpha _{k}\right\vert <n\right\}
}\lambda _{k}a_{k},\text{ \ }
\end{equation*}%
where%
\begin{equation}
\text{\ }\lambda _{k}=\frac{\Phi ^{1/2}\left( \alpha _{k}\right) }{%
V^{1/2}\left( \alpha _{k}\right) }\text{, \ \ }a_{k}=D_{2^{\left\vert \alpha
_{k}\right\vert +1}}-D_{2^{\left\vert \alpha _{k}\right\vert }}.  \label{101}
\end{equation}

Analogously to Theorem 1 if we apply Theorem W and (\ref{2aaa}) we conclude
that $F=\left( F_{n},n\in \mathbb{N}\right) \in H_{1}.$

By simple calculation we get that%
\begin{equation}
\widehat{F}(j)=\left\{ 
\begin{array}{ll}
\frac{\Phi ^{1/2}\left( \alpha _{k}\right) }{V^{1/2}\left( \alpha
_{k}\right) }, & \text{ if \thinspace \thinspace }j\in \left\{
2^{_{\left\vert \alpha _{k}\right\vert }},...,2^{_{\left\vert \alpha
_{k}\right\vert +1}}-1\right\} ,\text{ }k=0,1,... \\ 
0, & \text{\ if \thinspace \thinspace \thinspace }j\notin
\bigcup\limits_{k=0}^{\infty }\left\{ 2^{_{\left\vert \alpha _{k}\right\vert
}},...,2^{_{\left\vert \alpha _{k}\right\vert +1}}-1\right\} .\text{ }%
\end{array}%
\right.  \label{5aa}
\end{equation}

From (\ref{6aa}) and (\ref{5aa}) analogously to (\ref{6aaa}) we obtain%
\begin{equation*}
S_{\alpha _{k}}F=\sum_{\eta =0}^{k-1}\frac{\Phi ^{1/2}\left( \alpha _{\eta
}\right) }{V^{1/2}\left( \alpha _{\eta }\right) }\left( D_{2^{\left\vert
\alpha _{\eta }\right\vert +1}}-D_{2^{\left\vert \alpha _{\eta }\right\vert
}}\right) +\frac{\Phi ^{1/2}\left( \alpha _{k}\right) }{V^{1/2}\left( \alpha
_{k}\right) }w_{2^{\left\vert \alpha _{k}\right\vert }}D_{\alpha
_{k}-2^{\left\vert \alpha _{k}\right\vert }}.
\end{equation*}

By combining (\ref{22b}) and (\ref{2aaa}) we have 
\begin{equation*}
\left\Vert \frac{S_{\alpha _{k}}F}{\Phi \left( \alpha _{k}\right) }%
\right\Vert _{1}\geq \frac{\Phi ^{1/2}\left( \alpha _{k}\right) }{\Phi
\left( \alpha _{k}\right) V^{1/2}\left( \alpha _{k}\right) }\left\Vert
D_{_{\alpha _{k}-2^{\left\vert \alpha _{k}\right\vert }}}\right\Vert _{1}-%
\frac{1}{\Phi \left( \alpha _{k}\right) }\sum_{\eta =0}^{k-1}\frac{\Phi
^{1/2}\left( \alpha _{\eta }\right) }{V^{1/2}\left( \alpha _{\eta }\right) }
\end{equation*}%
\begin{equation*}
\geq \frac{V\left( \alpha _{k}-2^{\left\vert \alpha _{k}\right\vert }\right)
\Phi ^{1/2}\left( \alpha _{k}\right) }{8\Phi \left( \alpha _{k}\right)
V^{1/2}\left( \alpha _{k}\right) }-\frac{1}{\Phi \left( \alpha _{k}\right) }%
\sum_{\eta =0}^{\infty }\frac{\Phi ^{1/2}\left( \alpha _{\eta }\right) }{%
V^{1/2}\left( \alpha _{\eta }\right) }
\end{equation*}%
\begin{equation*}
\geq \frac{cV^{1/2}\left( \alpha _{k}\right) }{\Phi ^{1/2}\left( \alpha
_{k}\right) }\rightarrow \infty ,\text{ as }k\rightarrow \infty .
\end{equation*}

Theorem 2 is proved.
\end{proof}

\bigskip

\begin{proof}[\textbf{Proof of Theorem 3.}]
Let $0<p<1$ and $2^{k}<n\leq 2^{k+1}.$ By using Theorem 1 we see that%
\begin{equation}
\left\Vert S_{n}F-F\right\Vert _{H_{p}}\leq c_{p}\left\Vert
S_{n}F-S_{2^{k}}F\right\Vert _{H_{p}}+c_{p}\left\Vert
S_{2^{k}}F-F\right\Vert _{H_{p}}  \label{1000}
\end{equation}%
\begin{equation}
=c_{p}\left\Vert S_{n}\left( S_{2^{k}}F-F\right) \right\Vert
_{H_{p}}+c_{p}\left\Vert S_{2^{k}}F-F\right\Vert _{H_{p}}  \notag
\end{equation}%
\begin{equation*}
\leq c_{p}\left( 1+2^{d\left( n\right) \left( 1/p-1\right) }\right) \omega
_{H_{p}}\left( \frac{1}{2^{k}},F\right) \leq c_{p}2^{d\left( n\right) \left(
1/p-1\right) }\omega _{H_{p}}\left( \frac{1}{2^{k}},F\right) .
\end{equation*}

The proof of estimate (\ref{sn2}) is analogously to the proof of estimate (%
\ref{1000}). This completes the proof of theorem 3.
\end{proof}

\begin{proof}[\textbf{Proof of Theorem 4.}]
Let $0<p<1$, $F\in H_{p}$ and $\{m_{k}:k\geq 0\}$\ be a sequence of
nonnegative integers,\ satisfying condition (\ref{18a}). By using Theorem 3
we see that (\ref{con1}) holds.

Let proof second part of theorem 4. Under condition (\ref{dnk}), there
exists $\{\alpha _{k}:k\geq 0\}\subset \{m_{k}:k\geq 0\},$ such that 
\begin{equation}
\ 2^{d\left( \alpha _{k}\right) }\uparrow \infty ,\,\,\,\,\text{as}\ \
k\rightarrow \infty \text{, \ \ \ }2^{2\left( 1/p-1\right) d\left( \alpha
_{k}\right) }\leq 2^{\left( 1/p-1\right) d\left( \alpha _{k+1}\right) }.
\label{4.18}
\end{equation}%
We set 
\begin{equation*}
F_{n}=\sum_{\left\{ i:\text{ }\left\vert \alpha _{i}\right\vert <n\right\} }%
\frac{a_{i}}{2^{\left( 1/p-1\right) d\left( \alpha _{i}\right) }},
\end{equation*}%
where $a_{i}$ is defined by (\ref{100}). Since $a_{i}$ is $p$-atom if we
apply Theorem W and (\ref{4.18}) we conclude that $\ F\in H_{p}.$ On the
other hand%
\begin{eqnarray}
&&F-S_{2^{n}}F  \label{20} \\
&=&\left( F^{\left( 1\right) }-S_{2^{n}}F^{\left( 1\right) },...,F^{\left(
n\right) }-S_{2^{n}}F^{\left( n\right) },...,F^{\left( n+k\right)
}-S_{2^{n}}F^{\left( n+k\right) }\right)  \notag \\
&=&\left( 0,...,0,F^{\left( n+1\right) }-F^{\left( n\right) },...,F^{\left(
n+k\right) }-F^{\left( n\right) },...\right)  \notag \\
&=&\left( 0,...,0,\underset{i=n}{\overset{n+k}{\sum }}\frac{a_{i}}{2^{\left(
1/p-1\right) d\left( \alpha _{i}\right) }},...\right) ,\text{ \ }k\in 
\mathbb{N}_{+}  \notag
\end{eqnarray}%
is martingale. By combining (\ref{4.18}) and Theorem W we get that%
\begin{equation}
\omega _{H_{p}}(\frac{1}{2^{\left\vert \alpha _{k}\right\vert }},F)\leq
\sum\limits_{i=k}^{\infty }\frac{1}{2^{\left( 1/p-1\right) d\left( \alpha
_{i}\right) }}=O\left( \frac{1}{2^{\left( 1/p-1\right) d\left( \alpha
_{k}\right) }}\right) ,\text{ \ as \ }n\rightarrow \infty .  \label{4.21}
\end{equation}

It is easy to show that%
\begin{equation}
\widehat{F}(j)=\left\{ 
\begin{array}{ll}
2^{\left( 1/p-1\right) \left\langle \alpha _{k}\right\rangle }, & \text{ if
\thinspace \thinspace }j\in \left\{ 2^{_{\left\vert \alpha _{k}\right\vert
}},...,2^{_{\left\vert \alpha _{k}\right\vert +1}}-1\right\} ,\text{ }%
k=0,1,... \\ 
0, & \text{\ if \thinspace \thinspace \thinspace }j\notin
\bigcup\limits_{k=0}^{\infty }\left\{ 2^{_{\left\vert \alpha _{k}\right\vert
}},...,2^{_{\left\vert \alpha _{k}\right\vert +1}}-1\right\} .\text{ }%
\end{array}%
\right.  \label{4.22}
\end{equation}

Analogously to (\ref{77}) we can write that%
\begin{equation*}
\left\vert D_{\alpha _{k}}\right\vert \geq 2^{\left\langle \alpha
_{k}\right\rangle },\text{ \ \ \ \ for \ \ \ }I_{\left\langle \alpha
_{k}\right\rangle }\backslash I_{\left\langle \alpha _{k}\right\rangle +1}.
\end{equation*}

Since 
\begin{equation*}
\Vert D_{\alpha _{k}}\Vert _{L_{p,\infty }}\geq 2^{\left\langle \alpha
_{k}\right\rangle }\mu \left\{ x\in I_{\left\langle \alpha _{k}\right\rangle
}\backslash I_{\left\langle \alpha _{k}\right\rangle +1}:\text{ }\left\vert
D_{\alpha _{k}}\right\vert \geq 2^{\left\langle \alpha _{k}\right\rangle
}\right\} ^{1/p}
\end{equation*}%
\begin{equation*}
\geq 2^{\left\langle \alpha _{k}\right\rangle }\left( \mu \left\{
I_{\left\langle \alpha _{k}\right\rangle }\backslash I_{\left\langle \alpha
_{k}\right\rangle +1}\right\} \right) ^{1/p}\geq 2^{\left\langle \alpha
_{k}\right\rangle \left( 1-1/p\right) }
\end{equation*}%
by using (\ref{4.22}) we have%
\begin{equation*}
\Vert S_{\alpha _{k}}F-F\Vert _{L_{p,\infty }}\geq \Vert 2^{\left(
1/p-1\right) \left\langle \alpha _{k}\right\rangle }\left( D_{2^{\left\vert
\alpha _{k}\right\vert +1}}-D_{\alpha _{k}}\right) \Vert _{L_{p,\infty }}
\end{equation*}%
\begin{equation*}
-\left\Vert \sum\limits_{i=k+1}^{\infty }2^{\left( 1/p-1\right) \left\langle
\alpha _{i}\right\rangle }\left( D_{2^{\left\vert \alpha _{i}\right\vert
+1}}-D_{2^{\left\vert \alpha _{i}\right\vert }}\right) \right\Vert
_{L_{p,\infty }}
\end{equation*}%
\begin{equation*}
=2^{\left( 1/p-1\right) \left\langle \alpha _{k}\right\rangle }\Vert
D_{\alpha _{k}}\Vert _{L_{p,\infty }}-2^{\left( 1/p-1\right) \left\langle
\alpha _{k}\right\rangle }\Vert D_{2^{\left\vert \alpha _{k}\right\vert
+1}}\Vert _{L_{p,\infty }}
\end{equation*}%
\begin{equation*}
-\sum_{i\geq k+1}\frac{\Vert 2^{\left( 1/p-1\right) \left\vert \alpha
_{_{i}}\right\vert }\left( D_{2^{\left\vert \alpha _{_{i}}\right\vert
+1}}-D_{2^{\left\vert \alpha _{_{i}}\right\vert }}\right) \Vert
_{L_{p,\infty }}}{2^{\left( 1/p-1\right) d\left( \alpha _{i}\right) }}
\end{equation*}%
\begin{equation*}
\geq c-\frac{1}{2^{\left( 1/p-1\right) d\left( \alpha _{k}\right) }}%
-\sum_{i\geq k+1}\frac{1}{2^{\left( 1/p-1\right) d\left( \alpha _{i}\right) }%
}\geq c-\frac{2}{2^{\left( 1/p-1\right) d\left( \alpha _{k}\right) }}.
\end{equation*}

This completes the proof of Theorem 4.
\end{proof}

\begin{proof}[\textbf{Proof of Theorem 5.}]
\textbf{\ } Let $F\in H_{1}$ and $\{m_{k}:k\geq 0\}$\ be a sequence of
nonnegative integers,\ satisfying condition (\ref{cond1}). By using Theorem
3 we see that (\ref{cond1a}) holds.

Let proof second part of theorem 5. Under conditions of second part of
theorem 5, there exists sequence $\{\alpha _{k}:k\geq 0\}\subset
\{m_{k}:k\geq 0\}$ such that%
\begin{equation}
V(\alpha _{k})\uparrow \infty ,\,\,\,\,k\rightarrow \infty \text{ \ \ and\ \
\ \ }V^{2}(\alpha _{k})\leq V(\alpha _{k+1}).  \label{4.7}
\end{equation}

We set 
\begin{equation*}
F_{n}=\sum_{\left\{ i:\text{ }\left\vert \alpha _{i}\right\vert <n\right\} }%
\frac{a_{i}}{V(\alpha _{i})},
\end{equation*}%
where $a_{i}$ is defined by (\ref{101}). Since $a_{i}$ is a 1-atom if we
apply Theorem W and (\ref{4.7}) we conclude that $F=\left( F_{n},n\in 
\mathbb{N}\right) \in H_{1}.$

Analogously to (\ref{20}), (\ref{4.7}) and Theorem W we can show that%
\begin{equation*}
F-S_{2^{n}}F=\left( 0,...,0,\underset{i=n}{\overset{n+k}{\sum }}\frac{a_{i}}{%
V(\alpha _{i})},...\right) ,\text{ \ }k\in \mathbb{N}_{+}
\end{equation*}%
is martingale and%
\begin{equation}
\left\Vert F-S_{2^{n}}F\right\Vert _{H_{1}}\leq \sum\limits_{i=n+1}^{\infty }%
\frac{1}{V(\alpha _{i})}=O\left( \frac{1}{V(\alpha _{n})}\right) \text{\ \ \
as \ \ }n\rightarrow \infty .  \label{4.12}
\end{equation}

It is easy to show that%
\begin{equation}
\widehat{F}(j)=\left\{ 
\begin{array}{ll}
\frac{1}{V(\alpha _{k})}, & \text{ if \thinspace \thinspace }j\in \left\{
2^{_{\left\vert \alpha _{k}\right\vert }},...,2^{_{\left\vert \alpha
_{k}\right\vert +1}}-1\right\} ,\text{ }k=0,1,... \\ 
0, & \text{\ if \thinspace \thinspace \thinspace }j\notin
\bigcup\limits_{k=0}^{\infty }\left\{ 2^{_{\left\vert \alpha _{k}\right\vert
}},...,2^{_{\left\vert \alpha _{k}\right\vert +1}}-1\right\} .\text{ }%
\end{array}%
\right.  \label{4.13}
\end{equation}

By using (\ref{4.13}) we have 
\begin{equation*}
\Vert F-S_{\alpha _{k}}F\Vert _{1}\geq \Vert \frac{D_{2^{\left\vert \alpha
_{k}\right\vert +1}}-D_{\alpha _{k}}}{V(\alpha _{k})}+\overset{\infty }{%
\sum_{i=k+1}}\frac{D_{2^{\left\vert \alpha _{i}\right\vert
+1}}-D_{2^{\left\vert \alpha _{i}\right\vert }}}{V(\alpha _{i})}\Vert _{1}
\end{equation*}%
\begin{equation*}
\geq \frac{\Vert D_{\alpha _{k}}\Vert _{1}}{V(\alpha _{k})}-\frac{\Vert
D_{2^{\left\vert \alpha _{k}\right\vert +1}}\Vert _{1}}{V(\alpha _{k})}%
-\Vert \overset{\infty }{\sum_{i=k+1}}\frac{D_{2^{\left\vert \alpha
_{i}\right\vert +1}}-D_{2^{\left\vert \alpha _{i}\right\vert }}}{V(\alpha
_{i})}\Vert _{1}
\end{equation*}%
\begin{equation*}
\geq \frac{1}{8}-\frac{1}{V(\alpha _{k})}-\overset{\infty }{\sum_{i=k+1}}%
\frac{1}{V(\alpha _{i})}\geq \frac{1}{8}-\frac{2}{V(\alpha _{k})}.
\end{equation*}

This completes the proof of Theorem 5.
\end{proof}

\textbf{Acknowledgment: }The author would like to thank the referee for
helpful suggestions.

\end{document}